# Local polynomial regression on unknown manifolds

Peter J. Bickel[1] and Bo Li[2]

*University of California, Berkeley and Tsinghua University*

**Abstract:** We reveal the phenomenon that "naive" multivariate local polynomial regression can adapt to local smooth lower dimensional structure in the sense that it achieves the optimal convergence rate for nonparametric estimation of regression functions belonging to a Sobolev space when the predictor variables live on or close to a lower dimensional manifold.

## 1. Introduction

It is well known that worst case analysis of multivariate nonparametric regression procedures shows that performance deteriorates sharply as dimension increases. This is sometimes refered to as the curse of dimensionality. In particular, as initially demonstrated by [19, 20], if the regression function, $m(x)$, belongs to a Sobolev space with smoothness $p$, there is no nonparametric estimator that can achieve a faster convergence rate than $n^{-\frac{p}{2p+D}}$, where $D$ is the dimensionality of the predictor vector $X$.

On the other hand, there has recently been a surge in research on identifying intrinsic low dimensional structure from a seemingly high dimensional source, see [1, 5, 15, 21] for instance. In these settings, it is assumed that the observed high-dimensional data are lying on a low dimensional smooth manifold. Examples of this situation are given in all of these papers — see also [14]. If we can estimate the manifold, we can expect that we should be able to construct procedures which perform as well as if we know the structure. Even if the low dimensional structure obtains only in a neighborhood of a point, estimation at that point should be governed by actual rather than ostensible dimension. In this paper, we shall study this situation in the context of nonparametric regression, assuming the predictor vector has a lower dimensional smooth structure. We shall demonstrate the somewhat surprising phenomenon, suggested by Bickel in his 2004 Rietz lecture, that the procedures used with the expectation that the ostensible dimension $D$ is correct will, with appropriate adaptation not involving manifold estimation, achieve the optimal rate for manifold dimension $d$.

Bickel conjectured in his 2004 Rietz lecture that, in predicting $Y$ from $X$ on the basis of a training sample, one could automatically adapt to the possibility that the apparently high dimensional $X$ that one observed, in fact, lived on a much smaller dimensional manifold and that the regression function was smooth on that manifold. The degree of adaptation here means that the worst case analyses for prediction are governed by smoothness of the function on the manifold and not on the space in which $X$ ostensibly dwells, and that purely data dependent procedures can be constructed which achieve the lower bounds in all cases.

In this paper, we make this statement precise with local polynomial regression. Local polynomial regression has been shown to be a useful nonparametric technique in various local modelling, see [8, 9]. We shall sketch in Section 2 that local linear regression achieves this phenomenon for local smoothness $p = 2$, and will also argue that our procedure attains the global IMSE if global smoothness is assumed. We shall also sketch how polynomial regression can achieve the appropriate higher rate if more smoothness is assumed.

A critical issue that needs to be faced is regularization since the correct choice of bandwidth will depend on the unknown local dimension $d(x)$. Equivalently, we need to adapt to $d(x)$. We apply local generalized cross validation, with the help of an estimate of $d(x)$ due to [14]. We discuss this issue in Section 3. Finally we give some simulations in Section 4.

A closely related technical report, [2] came to our attention while this paper was in preparation. Binev et al consider in a very general way, the construction of non-parametric estimation of regression where the predictor variables are distributed according to a fixed completely unknown distribution. In particular, although they did not consider this possibility, their method covers the case where the distribution of the predictor variables is concentrated on a manifold. However, their method is, for the moment, restricted to smoothness $p \leq 1$ and their criterion of performance is the integral of pointwise mean square error with respect to the underlying distribution of the variables. Their approach is based on a tree construction which implicitly estimates the underlying measure as well as the regression. Our discussion is considerably more restrictive by applying only to predictors taking values in a low dimensional manifold but more general in discussing estimation of the regression function at a point. Binev et al promise a further paper where functions of general Lipschitz order are considered.

Our point in this paper is mainly a philosophical one. We can unwittingly take advantage of low dimensional structure without knowing it. We do not give careful minimax arguments, but rather, partly out of laziness, employ the semi heuristic calculations present in much of the smoothing literature.

Here is our setup. Let $(X_i, Y_i), (i = 1, 2, \ldots, n)$ be i.i.d $\Re^{D+1}$ valued random vectors, where $X$ is a $D$-dimensional predictor vector, $Y$ is the corresponding univariate response variable. We aim to estimate the conditional mean $m_0(x) = E(Y|X = x)$ nonparametrically. Our crucial assumption is the existence of a local *chart*, i.e., each small patch of $\mathcal{X}$ (a neighborhood around $x$) is isomorphic to a ball in a $d$-dimensional Euclidean space, where $d = d(x) \leq D$ may vary with $x$. Since we fix our working point $x$, we will use $d$ for the sake of simplicity. The same rule applies to other notations which may also depend on $x$.) More precisely, let $\mathcal{B}_{z,r}^d$ denote the ball in $\Re^d$, centered at $z$ with radius $r$. A similar definition applies to $\mathcal{B}_{x,R}^D$. For small $R > 0$, we consider the neighborhood of $x$, $\mathcal{X}_x := \mathcal{B}_{x,R}^D \cap \mathcal{X}$ within $\mathcal{X}$. We suppose there is a continuously differentiable bijective map $\phi : \mathcal{B}_{0,r}^d \mapsto \mathcal{X}_x$. Under this assumption with $d < D$, the distribution of $X$ degenerates in the sense that it does not have positive density around $x$ with respect to Lebesgue measure on $\Re^D$. However, the induced measure $\mathbb{Q}$ on $\mathcal{B}_{0,r}^d$ defined below, can have a non-degenerate density with respect to Lebesgue measure on $\Re^d$. Let $\mathcal{S}$ be an open subset of $\mathcal{X}_x$, and $\phi^{-1}(\mathcal{S})$ be its preimage in $\mathcal{B}_{0,r}^{d(x)}$. Then $\mathbb{Q}(Z \in \phi^{-1}(\mathcal{S})) = \mathbb{P}(X \in \mathcal{S})$. We assume throughout that $\mathbb{Q}$ admits a continuous positive density function $f(\cdot)$. We proceed to our main result whose proof is given in the Appendix.

## 2. Local linear regression

[17] develop the general theory for multivariate local polynomial regression in the usual context, i.e., the predictor vector has a $D$ dimensional compact support in $\Re^D$. We shall modify their proof to show the "naive" (brute-force) multivariate local linear regression achieves the "oracle" convergence rate for the function $m(\phi(z))$ on $\mathcal{B}_{0,r}^d$.

Local linear regression estimates the population regression function by $\hat{\alpha}$, where $(\hat{\alpha}, \hat{\beta})$ minimize

$$\sum_{i=1}^{n} \left(Y_i - \alpha - \beta^T(X_i - x)\right)^2 K_h(X_i - x).$$

Here $K_h(\cdot)$ is a $D-$variate kernel function. For the sake of simplicity, we choose the same bandwidth $h$ for each coordinate. Let

$$X_x = \begin{bmatrix} 1 & (X_1 - x)^T \\ \vdots & \vdots \\ 1 & (X_n - x)^T \end{bmatrix}$$

and $W_x = diag\{K_h(X_1 - x), \ldots, K_h(X_n - x)\}$. Then the estimator of the regression function can be written as

$$\hat{m}(x, h) = e_1^T (X_x^T W_x X_x)^{-1} X_x^T W_x Y$$

where $e_1$ is the $(D+1) \times 1$ vector having 1 in the first entry and 0 elsewhere.

### 2.1. Decomposition of the conditional MSE

We enumerate the assumptions we need for establishing the main result. Let $M$ be a canonical finite positive constant,

(i) The kernel function $K(\cdot)$ is continuous and radially symmetric, hence bounded.
(ii) There exists an $\epsilon(0 < \epsilon < 1)$ such that the following asymptotic irrelevance conditions hold.

$$E\left[K^\gamma(\frac{X-x}{h})w(X)1\left(X \in \left(\mathcal{B}_{x,h^{1-\epsilon}}^D \cap \mathcal{X}\right)^c\right)\right] = o(h^{d+2})$$

for $\gamma = 1, 2$ and $|w(x)| \leq M(1 + |x|^2)$.
(iii) $v(x) = Var(Y|X=x) \leq M$.
(iv) The regression function $m(x)$ is twice differentiable, and $\|\frac{\partial^2 m}{\partial x_a x_b}\|_\infty \leq M$ for all $1 \leq a \leq b \leq D$ if $x = (x_1, \ldots, x_D)$.
(v) The density $f(\cdot)$ is continuously differentiable and strictly positive at 0 in $\mathcal{B}_{0,r}^d$.

Condition (ii) is satisfied if $K$ has exponential tails since if $V = \frac{X-x}{h}$, the conditions can be written as

$$E\left[K^\gamma(V)w(x+hV)1(V \in (\mathcal{B}_{0,h^{1-\epsilon}}^D)^c\right] = o(h^{d+2}).$$

**Theorem 2.1.** *Let $x$ be an interior point in $\mathcal{X}$. Then under assumptions (i)-(v), there exist some $J_1(x)$ and $J_2(x)$ such that*

$$E\{\hat{m}(x,h) - m(x)|X_1,\ldots,X_n\} = h^2 J_1(x)(1 + o_P(1)),$$

$$Var\{\hat{m}(x,h) - m(x)|X_1,\ldots,X_n\} = n^{-1}h^{-d} J_2(x)(1 + o_P(1)).$$

**Remark 1.** The predictor vector doesn't need to lie on a perfect smooth manifold. The same conclusion still holds as long as the predictor vector is "close" to a smooth manifold. Here "close" means the noise will not affect the first order of our asymptotics. That is, we think of $X_1,\ldots,X_n$ as being drawn from a probability distribution $P$ on $\Re^D$ concentrated on the set

$$\mathcal{X} = \{y : |\phi(u) - y| \le \epsilon_n \text{ for some } u \in \mathcal{B}_{0,r}^d\}$$

and $\epsilon_n \to 0$ with $n$. It is easy to see from our arguments below that if $\epsilon_n = o(h)$, then our results still hold.

**Remark 2.** When the point of interest $x$ is on the boundary of the support $\mathcal{X}$, we can show that the bias and variance have similar asymptotic expansions, following the Theorem 2.2 in [17]. But, given the extra complication of the embedding, the proof would be messier, and would not, we believe, add any insight. So we omit it.

## 2.2. Extensions

It's somewhat surprising but not hard to show that if we assume the regression function $m$ to be $p$ times differentiable with all partial derivatives of order $p$ bounded ($p \ge 2$, an integer), we can construct estimates $\hat{m}$ such that,

$$E\{\hat{m}(x,h) - m(x)|X_1,\ldots,X_n\} = h^p J_1(x)(1 + o_P(1)),$$

$$Var\{\hat{m}(x,h) - m(x)|X_1,\ldots,X_n\} = n^{-1}h^{-d} J_2(x)(1 + o_P(1))$$

yielding the usual rate of $n^{-\frac{2p}{2p+d}}$ for the conditional MSE of $\hat{m}(x,h)$ if $h$ is chosen optimal, $h = \lambda n^{-\frac{1}{2p+d}}$. This requires replacing local linear regression with local polynomial regression with a polynomial of order $p-1$. We do not need to estimate the manifold as we might expect since the rate at which the bias term goes to 0 is derived by first applying Taylor expansion with respect to the original predictor components, then obtaining the same rate in the lower dimensional space by a first order approximation of the manifold map. Essentially all we need is that, locally, the geodesic distance is roughly proportionate to the Euclidean distance.

## 3. Bandwidth selection

As usual this tells us, for $p = 2$, that we should use bandwidth $\lambda n^{-\frac{1}{4+d}}$ to achieve the best rate of $n^{-\frac{2}{4+d}}$. This requires knowledge of the local dimension as well as the usual difficult choice of $\lambda$. More generally, dropping the requirement that the bandwidth for all components be the same we need to estimate $d$ and choose the constants corresponding to each component in a simple data determined way.

There is an enormous literature on bandwidth selection. There are three main approaches: plug-in ([7, 16, 18], etc); the bootstrap ([3, 11, 12], etc) and cross validation ([6, 10, 22], etc). The first has always seemed logically inconsistent to us since it requires higher order smoothness of $m$ than is assumed and if this higher order smoothness holds we would not use linear regression but a higher order polynomial. See also the discussion of [23].

We propose to use a blockwise cross-validation procedure defined as follows. Let the data be $(X_i, Y_i), 1 \leq i \leq n$. We consider a block of data points $\{(X_j, Y_j) : j \in \mathcal{J}\}$, with $|\mathcal{J}| = n_1$. Assuming the covariates have been standardized, we choose the same bandwidth $h$ for all the points and all coordinates within the block. A leave-one-out cross validation with respect to the block while using the whole data set is defined as following. For each $j \in \mathcal{J}$, let $\hat{m}_{-j,h}(X_j)$ be the estimated regression function (evaluated at $X_j$) via local linear regression with the whole data set except $X_j$. In contrast to the usual leave-one-out cross-validation procedure, our modified leave-one-out cross-validation criterion is defined as $mCV(h) = \frac{1}{n_1} \sum_{j \in \mathcal{J}} (Y_j - \hat{m}_{-j,h}(X_j))^2$. Using a result from [23], it can be shown that

$$mCV(h) = \frac{1}{n_1} \sum_{j \in \mathcal{J}} \frac{(Y_j - \hat{m}_h(X_j))^2}{(1 - S_h(j,j))^2}$$

where $S_h(j,j)$ is the diagonal element of the smoothing matrix $S_h$. We adopt the GCV idea proposed by [4] and replace the $S_h(j,j)$ by their average $atr_{\mathcal{J}}(S_h) = \frac{1}{n_1} \sum_{j \in \mathcal{J}} S_h(j,j)$. Thereby our modified generalized cross-validation criterion is,

$$mGCV(h) = \frac{1}{n_1} \sum_{j \in \mathcal{J}} \frac{(Y_j - \hat{m}_h(X_j))^2}{(1 - atr_{\mathcal{J}}(S_h))^2}.$$

The bandwidth $h$ is chosen to minimize this criterion function.

We give some heuristics for the justifying the (blockwise homoscedastic) mGCV. In a manner analogous to [23], we can show

$$S_h(j,j) = e_1^T (X_x^T W_x X_x)^{-1} e_1 K_h(0)|_{x=X_j}.$$

In view of (A.2) in the Appendix, we see $S_h(j,j) = n^{-1} h^{-d} K(0)(A_1(X_j) + o_p(1))$. Thus as $n^{-1} h^{-d} \to 0$,

$$atr_{\mathcal{J}}(S_h) = n^{-1} h^{-d} K(0)(n_1^{-1} \sum_{j \in \mathcal{J}} A_1(X_j) + o_p(1))$$
$$= O_p(n^{-1} h^{-d}) = o_p(1).$$

Then, as is discussed in [22], using the approximation $(1-x)^{-2} \approx 1 + 2x$ for small $x$, we can rewrite $mGCV(h)$ as

$$mGCV(h) = \frac{1}{n_1} \sum_{j \in \mathcal{J}} (Y_j - \hat{m}_h(X_j))^2 + \frac{2}{n_1} tr_{\mathcal{J}}(S_h) \frac{1}{n_1} \sum_{j \in \mathcal{J}} (Y_j - \hat{m}_h(X_j))^2.$$

Now regarding $\frac{1}{n_1} \sum_{j \in \mathcal{J}} (Y_j - \hat{m}_h(X_j))^2$ in the second term as an estimator of the constant variance for the focused block, the mGCV is approximately the same as the $C_p$ criterion, which is an estimator of the prediction error up to a constant.

In practice, we first use [14]'s approach to estimate the local dimension $d$, which yields a consistent estimate $\hat{d}$ of $d$. Based on the estimated intrinsic dimensionality $\hat{d}$, a set of candidate bandwidths $\mathcal{CB} = \{\lambda_1 n^{-\frac{1}{d+4}}, \ldots, \lambda_B n^{-\frac{1}{d+4}}\}$ ($\lambda_1 < \cdots < \lambda_B$) are chosen . We pick the one minimizing the $mGCV(h)$ function.

## 4. Numerical experiments

The data generating process is as following. The predictor vector $X = (X_{(1)}, X_{(2)}, X_{(3)})$, where $X_{(1)}$ will be sampled from a standard normal distribution, $X_{(2)} = X_{(1)}^3 + sin(X_{(1)}) - 1$, and $X_{(3)} = \log(X_{(1)}^2 + 1) - X_{(1)}$. The regression function $m(x) = m(x_{(1)}, x_{(2)}, x_{(3)}) = cos(x^{(1)}) + x_{(2)} - x_{(3)}^2$. The response variable $Y$ is generated via the mechanism $Y = m(X) + \varepsilon$, where $\varepsilon$ has a standard normal distribution. By definition, the 3-dimensional regression function $m(x)$ is essentially a 1-dimensional function of $x_{(1)}$. $n = 200$ samples are drawn. The predictors are standardized before estimation. We estimate the regression function $m(x)$ by both the "oracle" univariate local linear (ull) regression with a single predictor $X_{(1)}$ and our blind 3-variate local linear regression with all predictors $X_{(1)}, X_{(2)}, X_{(3)}$.

We focus on the middle block with 100 data points, with the number of neighbor parameter $k$, needed for Levina and Bickel's estimate, set to be 15. The intrinsic dimension estimator is $\hat{d} = 1.023$, which is close to the true dimension, $d = 1$. We use the Epanechnikov kernel in our simulation. Our proposed modified GCV procedure is applied to both the ull and mll procedures. The estimation results are displayed in Figure 1. The $x-axis$ is the standardized $X_{(1)}$. From the right panel, we see the blind mll indeed performs almost as well as the "oracle" ull.

Next, we allow the predictor vector to only lie close to a manifold. Specifically, we sample $X_{(1)} = X'_{(1)} + \epsilon'_1, X_{(2)} = X'^3_{(1)} + sin(X'_{(1)}) - 1 + \epsilon'_2, X_{(3)} = \log(X'^2_{(1)} + 1) - X'_{(1)} + \epsilon'_3$, where $X'_{(1)}$ is sampled from a standard normal distribution, and $\epsilon'_1, \epsilon'_2$ and $\epsilon'_3$ are sampled from $\mathcal{N}(0, \sigma'^2)$. The noise scale is hence governed by $\sigma'$. In our experiment, $\sigma'$ is set to be $0.02, 0.04, \ldots, 0.18, 0.20$ respectively. The predictor vector samples are visualized in the left panel of Figure 2 with $\sigma' = 0.20$. In the maximum noise scale case, the pattern of the predictor vector is somewhat vague. Again, a blind "mll" estimation is done with respect to new data generated in the aforementioned way. We plot the MSEs associated with different noise scales in the right panel of Figure 2. The moderate noise scales we've considered indeed don't have a significant influence on the performance of the "mll" estimator in terms of MSE.

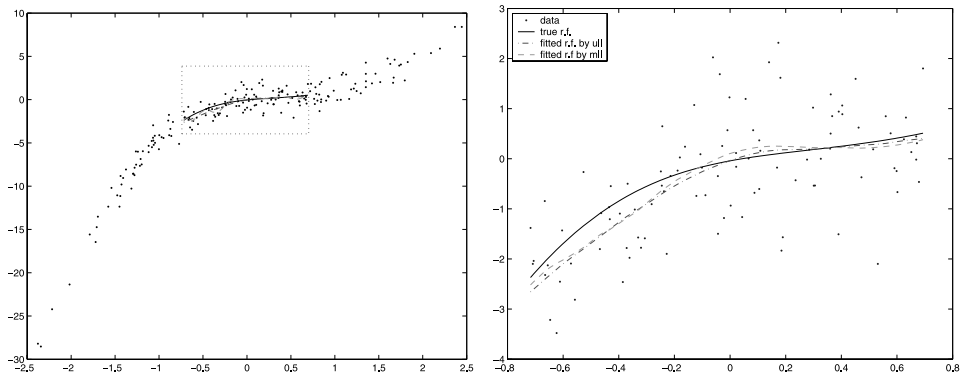

FIG 1. *The case with perfect embedding. The left panel shows the complete data and fitting of the middle block by both univariate local linear (ull) regression and multivariate local linear (mll) regression with bandwidths chosen via our modified GCV. The focused block is amplified in the right panel.*

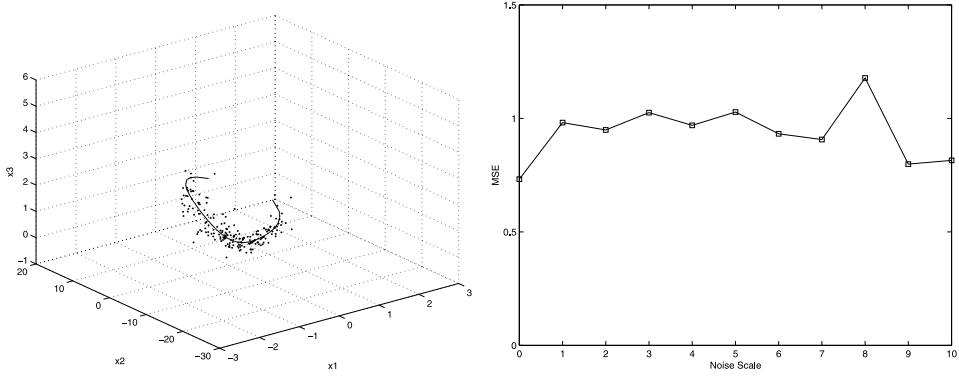

FIG 2. *The case with "imperfect" embedding. The left panel shows the predictor vector in a 3-D fashion with the noise scale $\sigma' = 0.2$. The right panel gives the MSEs with respect to increasing noise scales.*

## Appendix

*Proof of Theorem 2.1.* Using the notation of [17], $\mathcal{H}_m(x)$ is the $D \times D$ Hessian matrix of $m(x)$ at $x$, and

$$Q_m(x) = [(X_1 - x)^T \mathcal{H}_m(x)(X_1 - x), \cdots, (X_n - x)^T \mathcal{H}_m(x)(X_n - x)]^T.$$

Ruppert and Wand have obtained the bias term.

(A.1) $$\begin{aligned} &E(\hat{m}(x,h) - m(x)|X_1, \cdots, X_n) \\ &= \frac{1}{2} e_1^T (X_x^T W_x X_x)^{-1} X_x^T W_x \{Q_m(x) + R_m(x)\} \end{aligned}$$

where if $|\cdot|$ denotes Euclidean norm, $|R_m(x)|$ is of lower order than $|Q_m(x)|$. Also we have

$$\begin{aligned} &n^{-1} X_x^T W_x X_x \\ &= \begin{bmatrix} n^{-1} \sum_{i=1}^n K_h(X_i - x) & n^{-1} \sum_{i=1}^n K_h(X_i - x)(X_i - x)^T \\ n^{-1} \sum_{i=1}^n K_h(X_i - x)(X_i - x) & n^{-1} \sum_{i=1}^n K_h(X_i - x)(X_i - x)(X_i - x)^T \end{bmatrix}. \end{aligned}$$

The difference in our context lies in the following asymptotics.

$$\begin{aligned} EK_h(X_i - x) &= E\big[K_h(X_i - x)\mathbf{1}\big(X_i \in \mathcal{B}_{x,h^{1-\epsilon}}^D \cap \mathcal{X}\big)\big] \\ &\quad + E\big[K_h(X_i - x)\mathbf{1}\big(X_i \in \big(\mathcal{B}_{x,h^{1-\epsilon}}^D \cap \mathcal{X}\big)^c\big)\big] \\ &\stackrel{(ii)}{=} h^{-D}\Big(\int_{N_{0,h^{1-\epsilon}}^d} K\Big(\frac{\phi(z') - \phi(0)}{h}\Big) f(z') dz' + o_P(h^d)\Big) \\ &= h^{d-D}\Big(f(0) \int_{\Re^d} K(\nabla\phi(0)u) du + o_P(1)\Big) \\ &= h^{d-D}\big(A_1(x) + o_P(1)\big). \end{aligned}$$

Thus, by the LLN, we have

$$n^{-1} \sum_{i=1}^n K_h(X_i - x) = h^{d-D}\big(A_1(x) + o_P(1)\big).$$

Similarly, there exist some $A_2(x)$ and $A_3(x)$ such that

$$n^{-1} \sum_{i=1}^{n} K_h(X_i - x)(X_i - x) = h^{2+d-D}(A_2(x) + o_P(1))$$

and

$$n^{-1} \sum_{i=1}^{n} K_h(X_i - x)(X_i - x)(X_i - x)^T = h^{2+d-D}(A_3(x) + o_P(1))$$

where we used assumption (i) to remove the term of order $h^{1+d-D}$ in deriving the asymptotic behavior of $n^{-1} \sum_{i=1}^{n} K_h(X_i - x)(X_i - x)$. Invoking Woodbury's formula, as in the proof of Lemma 5.1 in [13], leads us to

(A.2) $$\left(n^{-1} X_x^T W_x X_x\right)^{-1} = h^{D-d} \begin{bmatrix} A_1(x)^{-1} + o_P(1) & O_P(1) \\ O_P(1) & h^{-2} O_p(1) \end{bmatrix}$$

On the other hand,

$$n^{-1} X_x W_x Q_m(x) = \begin{bmatrix} n^{-1} \sum_{i=1}^{n} K_h(X_i - x)(X_i - x)^T \mathcal{H}_m(x)(X_i - x) \\ n^{-1} \sum_{i=1}^{n} \{K_h(X_i - x)(X_i - x)^T \mathcal{H}_m(x)(X_i - x)\}(X_i - x) \end{bmatrix}.$$

In a similar fashion, we can deduce that for some $B_1(x), B_2(x)$,

$$n^{-1} \sum_{i=1}^{n} K_h(X_i - x)(X_i - x)^T \mathcal{H}_m(x)(X_i - x) = h^{2+d-D}(B_1(x) + o_P(1))$$

and

$$n^{-1} \sum_{i=1}^{n} \{K_h(X_i - x)(X_i - x)^T \mathcal{H}_m(x)(X_i - x)\}(X_i - x) = h^{3+d-D}(B_2(x) + o_P(1)).$$

We have

(A.3) $$n^{-1} X_x W_x Q_m(x) = h^{d-D} \begin{bmatrix} h^2(B_1(x) + o_P(1)) \\ h^3(B_2(x) + o_P(1)) \end{bmatrix}.$$

It follows from (A.1),(A.2) and (A.3) that the bias admits the following approximation.

(A.4) $$E(\hat{m}(x,h) - m(x)|X_1, \ldots, X_n) = h^2 A_1(x)^{-1} B_1(x) + o_P(h^2).$$

Next, we move to the variance term.

(A.5) $$\begin{aligned} &Var\{\hat{m}(x,h)|X_1, \ldots, X_n\} \\ &= e_1^T (X_x^T W_x X_x)^{-1} X_x^T W_x V W_x X_x (X_x^T W_x X_x)^{-1} e_1. \end{aligned}$$

The upper-left entry of $n^{-1} X_x^T W_x V W_x X_x$ is

$$n^{-1} \sum_{i=1}^{n} K_h(X_i - x)^2 v(X_i) = h^{d-2D} C_1(x)(1 + o_P(1)).$$

The upper-right block is

$$n^{-1} \sum_{i=1}^{n} K_h(X_i - x)^2 (X_i - x)^T v(X_i) = h^{1+d-2D} C_2(x)(1 + o_P(1))$$

and the lower-right block is

$$n^{-1} \sum_{i=1}^{n} K_h(X_i - x)^2 (X_i - x)(X_i - x)^T v(X_i) = h^{2+d-2D} C_3(x)(1 + o_P(1)).$$

In light of (A.2), we arrive at

(A.6) $\quad Var\{\hat{m}(x,h)|X_1,\ldots,X_n\} = n^{-1} h^{-d} A_1(x)^{-2} C_1(x)(1 + o_P(1)).$

The proof is complete. $\square$

**Acknowledgment.** We thank Ya'acov Ritov for insightful comments.